\documentclass[a4paper,11pt]{article} 
\usepackage{float}
\usepackage{booktabs} 
\usepackage{tikz-qtree}
\usetikzlibrary{arrows.meta}
\usepackage[a4paper, margin=.9in]{geometry}
\usepackage{amsmath,amsthm,amsfonts,amssymb,array}
\usepackage[square,sort,comma,numbers]{natbib}
\usepackage{subeqnarray}
\usepackage{blkarray}
\usepackage[utf8]{inputenc}
\usepackage[english]{babel}
\usepackage{mathrsfs}
\usepackage{colortbl}
\usepackage{graphicx,epsfig}
\newcounter{myeqno}
\usepackage{color}
\usepackage{tcolorbox}
\usepackage{ upgreek }
\usepackage{ mathrsfs }
\usepackage{relsize}
\usepackage{multirow}
\definecolor{shadecolor}{gray}{0.75}

\theoremstyle{definition}
\newtheorem{definition}{Definition}[section]

\newtheorem{proposition}{Proposition}[section]

\newtheorem{theorem}{Theorem}[section]
\newtheorem{example}{Example}[section]
\newtheorem{corollary}{Corollary}[section]

\newtheorem{lemma}{Lemma}[section]

\usepackage[misc,geometry]{ifsym}
\date{}
\usepackage[colorlinks=true]{hyperref}
\hypersetup{urlcolor=blue, citecolor=blue}
\newlength{\defbaselineskip}
\newcommand{\setlinespacing}[1]%
{\setlength{\baselineskip}{#1 \defbaselineskip}}

\begin{document}

		\title{\textbf{Resistance distance in $k$-coalescence of certain graphs}}

\author{Haritha T$^1$\footnote{harithathottungal96@gmail.com},  Chithra A. V$^1$\footnote{chithra@nitc.ac.in}
	\\ \\ \small 
 1 Department of Mathematics, National Institute of Technology Calicut,\\\small Calicut-673 601, Kerala, India\\ \small}

\maketitle	
\begin{abstract}
Any graph can be considered as a network of resistors, each of which has a resistance of $1 \Omega.$ The resistance distance $r_{ij}$ between a pair of vertices $i$ and $j$ in a graph is defined as the effective resistance between $i$ and $j$. This article deals with the resistance distance in the $k$-coalescence of complete graphs. We also present its results in connection with the Kemeny's constant, Kirchhoff index, additive degree-Kirchhoff index, multiplicative degree-Kirchhoff index and mixed degree-Kirchhoff index. Moreover, we obtain the resistance distance in the $k$-coalescence of a complete graph with particular graphs. As an application, we provide the resistance distance of certain graphs such as the vertex coalescence of a complete bipartite graph with a complete graph, a complete bipartite graph with a star graph, the windmill graph, pineapple graph, etc.
\end{abstract}
\section{Introduction}

Let $G_n= (V(G_n),E(G_n))$ be a simple connected undirected graph, consisting of $n$ vertices $\{v_1, v_2, \ldots, v_n\}$ and $m$ edges $\{e_1, e_2, \ldots, e_m\}$. The \textit{adjacency matrix} $A(G_n) = (a_{ij})$ of $G_n$ is defined such that $a_{ij} = 1$ if vertex $v_i$ is adjacent to vertex $v_j$, and it is zero otherwise. Denote the all-one entry matrix by $J_{n\times m}$, and the identity matrix by $I_{n}$.  The \textit{complete graph}, \textit{path}, \textit{cycle} and \textit{star graph} are denoted by $K_n$, $P_n$, $C_n$, and $K_{1,n}$ respectively, and $K_{n_1, n_2}$ is said to be the complete bipartite graph. Let $d_i$ denote the degree of a vertex $v_i$ in $G_n$. Note that the \textit{Laplacian matrix} $L(G_n)= D(G_n)-A(G_n)$, where $D(G_n)$ is the diagonal matrix of vertex degrees.
 The concept of \textit{resistance distance} is introduced by Klein and Randić in 1993 \cite{klein1993resistance}. The authors presented a new point of view, if we assign fixed resistances to each edge of a connected graph, then the resulting effective resistance between pairs of vertices corresponds to a graphical distance. For an $m\times n$ matrix, the matrix $P$ of order $n\times m$ is said to be a $\{1\}$-inverse of $M$ (denoted by $M^{(1)}$)  if $MPM = M$.  For any square matrix $N$, its group-inverse $N^{\#}$, refers to a distinct matrix $X$ that satisfies three conditions: $NXN= N$, $XNX= X$, and $NX= XN$. Clearly, the group inverse of $N$ is a $\{1\}$-inverse of $N$ \cite{ben2003generalized}.

 The standard method to compute the resistance distance $r_{ij}$ \cite{bapat2010graphs} between two vertices $v_i$ and $v_j$ is by using the $\{1\}$-inverse and group inverse of the Laplacian matrix $L=(l_{ij})$ of the underlying graph $G_n$ which is

 \begin{equation*}
     \begin{aligned}
         r_{ij} &= l^{(1)}_{ii}+ l^{(1)}_{jj}-l^{(1)}_{ij}-l^{(1)}_{ji}=l_{ii}^{\#}+l_{jj}^{\#}-2l_{ij}^{\#}.
     \end{aligned}
 \end{equation*}

The matrix $R(G_n)= (r_{ij})_{n\times n}$ is called the resistance distance matrix of $G_n$. The \textit{resistance distance energy} $RE(G_n)$ of $G_n$ is defined as the sum of the absolute
 values of the eigenvalues of $R(G_n).$

\par \textit{Kemeny's constant}  \cite{kirkland2012group} is essential in the theory of random walks, and it measures the average time it takes for a random walk to reach a vertex. It is defined as
 $$\kappa(G_n)= \frac{1}{4m}\sum_{v_i, v_j\in V(G_n)}d_{i}d_{j}r_{ij}.$$
 The \textit{Kirchhoff index} of $G_n$, also known as the total resistance of a network, represented as $\mathcal{K}f(G_n)$ \cite{klein1993resistance,bonchev1994molecular},
is defined as,
$$\mathcal{K}f(G_n)= \sum_{i<j} r_{ij}.$$
 The following are three graph parameters which are in terms of vertex degrees and resistance distance of a graph $G_n.$\\

The \textit{mixed degree-Kirchhoff index} of $G_n$ \cite{bianchi2016upper} is
$$\hat{R}(G_n)= \sum_{i<j}\left(\frac{d_i}{d_j}+\frac{d_j}{d_i}\right)r_{ij}.$$

 The \textit{multiplicative degree-Kirchhoff index} \cite{chen2007resistance} of $G_n$ is 
 $$R^{*}(G_n)= \sum_{i<j}d_id_jr_{ij}.$$
 
 The \textit{additive degree-Kirchhoff index} \cite{gutman2012degree} of $G_n$ is
 $$R^{+}(G_n)= \sum_{i<j}(d_i+d_j)r_{ij}.$$
 
 \par Suppose we have two graphs $G_n$ and $G'_{n'}$ with $v\in V(G_n)$ and $v'\in V(G'_{n'})$, then the coalescence $G_n \circ_{1} G'_{n'}$\cite{cvet} of $G_n$ and $G'_{n'}$ with respect to $v$ and $v'$ is formed by identifying $v$ and $v'$. Sudhir et al. introduced the concept of $k$-coalescence of two graphs in their work \cite{kcoal}, and it is defined as follows:
\begin{definition}
Let $G_n$ and $G'_{n'}$ be two connected graphs of orders $n$ and $n'$ and sizes $m$ and $m'$ respectively having an induced complete graph order $k$ with $n, n'\geq k$. Then the $k$-coalescence $G_n \circ_{k}G'_{n'}$ of $G_n$ and $G'_{n'}$ is the graph obtained
by identifying $k$ vertices on $^kC_{2}$ edges of induced $K_{k}$. The order and size of $G_n \circ_{k}G'_{n'}$ are  $n+n'-k$ and $m+m'-^kC_{2}$ respectively.
\end{definition}

\begin{figure}[htbp]
\centering
 \begin{tikzpicture}[scale=0.4,inner sep=1.2pt]
\draw (2,0) node(4) [circle,draw,fill] {}
      (16,0) node(10) [circle,draw,fill] {}
      (0,2) node(3) [circle,draw,fill] {}
      (4,2) node(1) [circle,draw,fill] {}
      (6,1) node(5) [circle,draw,fill] {}
      (6,3) node(6) [circle,draw,fill] {}
      (2,4) node(2) [circle,draw,fill] {}
      (16,4) node(8) [circle,draw,fill] {}
      (14,2) node(9) [circle,draw,fill] {}
      (18,2) node(7) [circle,draw,fill] {}
      (20,1) node(11) [circle,draw,fill] {}
      (20,3) node(12) [circle,draw,fill] {};
            
\draw [-] (1) to (2) to (3) to (4) to (1) to (5) to (6);
\draw[-] (6) to (1) to (3);
\draw[-] (2) to (4);
\draw[-] (12) to (7) to (10) to (9) to (8) to (7) to (11) ;
\draw[-] (9) to (7) to (12);
\draw[-] (8) to (10);

\end{tikzpicture}
 \caption{$K_4\circ K_3$\; \text{and}\; $K_4\circ K_{1,2}.$}
 \end{figure}
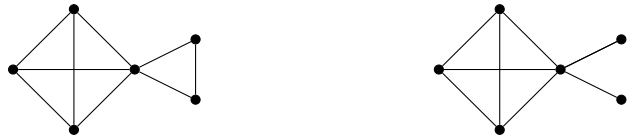

\par Resistance distance is significant in combinatorial matrix
theory \cite{bapat2010graphs, bapat2010resistance} and spectral graph theory \cite{bapat1999resistance,chen2007resistance,bapat2003simple,zhang2007resistance,subdivision}. For a survey of methods for finding resistance distance in graphs see \cite{evans2022algorithmic}. Graph operations have been widely used to analyse complex networks with properties abstracted from the real world. The formulas for resistance distance and Kirchhoff index pertaining to numerous graph classes and graph operations were presented in \cite{subdivision,liu2015resistance,wang2019resistance}. The aim of this work is to give the resistance distance in the $k$-coalescence of complete graphs and to provide parameters such as the Kemeny's constant, Kirchhoff index, additive degree-Kirchhoff index, multiplicative degree-Kirchhoff index and mixed degree-Kirchhoff index. Moreover, we obtain these results for some class of graphs.

\section{Preliminaries}
Through this section we present some useful lemmas and theorems.
\begin{lemma}\cite{zhou2014line}\label{P}
Let $C= \begin{bmatrix}C_0& C_1\\ C_2& C_3\end{bmatrix}$ be a nonsingular matrix. If $C_0$ and $C_3$ are nonsingular, then
\begin{equation*}
    \begin{aligned}
    C^{-1}&= \begin{bmatrix}(C_{0}-C_{1}C_{3}^{-1}C_{2})^{-1}& -C_{0}^{-1}C_{1}P^{-1}\\ -P^{-1}C_{2}C_{0}^{-1}& P^{-1}\end{bmatrix},
    \end{aligned}
\end{equation*}
where $P= C_{3}-C_{2}C_{0}^{-1}C_{1}.$

\end{lemma}

\begin{lemma}\cite{subdivision}\label{L}
Let $L= \begin{bmatrix}L_1& L_2\\ L_{2}^{T}& L_3\end{bmatrix}$ be the Laplacian matrix of a connected graph. If each column vector of $L_{2}^{T}$ is $-e$ (the all-ones column vector) or a zero vector, then $L^{(1)}= \begin{bmatrix}L_{1}^{-1}& 0\\ 0& S^{\#}\end{bmatrix}$, where $S= L_{3}-L_{2}^{T}L_{1}^{-1}L_{2}.$
\end{lemma}
\begin{lemma}\cite{subdivision}\label{dag}
Let $L$ be the Laplacian matrix of a graph $G_n$. For any $a >0$, we have 
$(L +aI-\frac{a}{n}J_n)^{\#}= (L +aI)^{-1}-\frac{1}{an}J_n.$
\end{lemma}
\begin{lemma}\cite{jahfar2020central}\label{rs}
For any real numbers $r,s>0$,
$$(rI_{n}-sJ_{n})^{-1}= \frac{1}{r}I_{n}+\frac{s}{r(r-ns)}J_{n}.$$
\end{lemma}

\section{Main results}
This section provides the resistance distance in $k$-coalescence of certain graphs and discusses some of its graph parameters. The graph $K_{p_1}\circ_{k}K_{p_2}$ consists of $t= p_1+p_2-k$ vertices.
\begin{theorem}\label{rd}
For $p_1, p_2\geq k$, let $T$ be the collection of vertices in $K_{p_1}\circ_{k}K_{p_2}$, which are the identified vertices of some vertices in $K_{p_1}$ and $K_{p_2}$. Then,
\begin{itemize}
    \item [(i)] for any $v_i, v_j\in T$,
    $r_{ij}= \frac{2}{t}.$
    \item [(ii)] for $v_i\in T, v_j\in V(K_{p_1}\setminus T)$,
    $r_{ij}= \frac{(k+1)(p_2-k)+2p_1k}{kp_1t}.$
    \item [(iii)] for $v_i\in T, v_j\in V(K_{p_2}\setminus T)$,
    $r_{ij}= \frac{(k+1)(p_1-k)+2p_2k}{kp_2t}.$
    \item [(iv)] for any $v_i, v_j\in V(K_{p_1}\setminus T), r_{ij}= \frac{2}{p_1}.$
    \item [(v)] for $v_i\in V(K_{p_1}\setminus T), v_j\in V(K_{p_2}\setminus T), r_{ij}= \frac{(p_1+p_2)(k+1)}{kp_1p_2}.$
    \item [(vi)] for any $v_i, v_j\in V(K_{p_2}\setminus T), r_{ij}= \frac{2}{p_2}.$
\end{itemize}

\end{theorem}
\begin{proof}
The Laplacian matrix of $K_{p_1}\circ_k K_{p_2}$ is given by,
$$L(K_{p_1}\circ_k K_{p_2})= \begin{bmatrix}tI_k-J_k& -J_{k\times p_1-k}& -J_{k\times p_2-k}\\-J_{k\times p_1-k}^T& p_1 I_{p_{1}-k}-J_{p_1-k}& 0\\-J_{k\times p_2-k}^T& 0& p_2 I_{p_{2}-k}-J_{p_2-k}\end{bmatrix}.$$
Let $L_1= \begin{bmatrix}tI_k-J_k& -J_{k\times p_1-k}\\-J_{k\times p_1-k}^T& p_1 I_{p_1-k}-J_{p_1-k}\end{bmatrix}, L_2= \begin{bmatrix}-J_{k\times p_2-k}\\0\end{bmatrix}$, and $L_3= p_2I-J_{p_2-k}.$\\
Then by Lemmas \ref{P} and \ref{rs} we get,\\
$$L_{1}^{-1}= \begin{bmatrix}\frac{1}{t}(I_k+\frac{p_1}{k(p_2-k)}J_k)& \frac{1}{k(p_2-k)}J_{k\times p_1-k}\\\frac{1}{k(p_2-k)}J_{k\times p_1-k}^T& \frac{1}{p_1}I_{p_1-k}+\frac{t}{p_1k(p_2-k)}J_{p_1-k}\end{bmatrix}.$$
Consider 
\begin{equation*}
    \begin{aligned}
    S&= L_3-L_{2}^{T}L_{1}^{-1}L_2\\
    &= p_2I_{p_2-k}-\frac{p_2}{p_2-k}J_{p_2-k}
    \end{aligned}
\end{equation*}

then, $S^\#= \frac{1}{p_2}(I_{p_2-k}-\frac{1}{p_2-k}J_{p_2-k}).$\\
From Lemma \ref{L},
$$L^{(1)}(K_{p_1}\circ_{k}K_{p_2})= \begin{bmatrix}\frac{1}{t}(I_k+\frac{p_1}{k(p_2-k)}J_k)& \frac{1}{k(p_2-k)}J_{k\times p_1-k}& 0\\\frac{1}{k(p_2-k)}J_{k\times p_1-k}^T& \frac{1}{p_1}I_{p_1-k}+\frac{t}{p_1k(p_2-k)}J_{p_1-k}& 0\\0& 0& \frac{1}{p_2}(I_{p_2-k}-\frac{1}{p_2-k}J_{p_2-k})\end{bmatrix}.$$
For any $v_i, v_j\in T$,
\begin{equation*}
    \begin{aligned}
    r_{ij}&= \frac{2}{t}\left(1+\frac{p_1}{k(p_2-k)}\right)-\frac{2p_1}{kt(p_2-k)}\\
    &= \frac{2}{t}.
    \end{aligned}
\end{equation*}\\

For $v_i\in T, v_j\in V(K_{p_1}\setminus T)$,
\begin{equation*}
    \begin{aligned}
    r_{ij}&= \frac{1}{t}\left(1+\frac{p_1}{k(p_2-k)}\right)-\frac{1}{p_1}\left(1+\frac{t}{k(p_2-k)}\right)-\frac{2}{k(p_2-k)}\\
    &= \frac{(k+1)(p_2-k)+2p_1k}{kp_1t}.
    \end{aligned}
\end{equation*}\\

For $v_i\in T, v_j\in V(K_{p_2}\setminus T)$,
\begin{equation*}
    \begin{aligned}
    r_{ij}&= \frac{1}{t}+\frac{p_1}{k(p_2-k)t}+\frac{1}{p_2}-\frac{1}{p_2(p_2-k)}\\
    &= \frac{(k+1)(p_1-k)+2p_2k}{kp_2kt}.
    \end{aligned}
\end{equation*}

 For any $v_i, v_j\in V(K_{p_1}\setminus T), r_{ij}= \frac{2}{p_1}.$\\
 
 For $v_i\in V(K_{p_1}\setminus T), v_j\in V(K_{p_2}\setminus T)$,
  \begin{equation*}
      \begin{aligned}
      r_{ij}&= \frac{1}{p_1}\left(1+\frac{t}{k(p_2-k)}\right)+\frac{1}{p_2}\left(1-\frac{1}{p_2-k}\right)\\
      &= \frac{(p_1+p_2)(k+1)}{kp_1p_2}.
      \end{aligned}
  \end{equation*}
For any $v_i,\; v_j\in V(K_{p_2}\setminus T), r_{ij}= \frac{2}{p_2}.$\\
Therefore, the resistance distance matrix of $K_{p_1}\circ_{k}K_{p_2}$ is 
$$R(K_{p_1}\circ_{k}K_{p_2})= \begin{bmatrix}\frac{2}{t}(J_k-I_k)& \frac{(k+1)(p_2-k)+2p_1k}{kp_1t}& \frac{(k+1)(p_1-k)+2p_2k}{p_2kt}J_{k\times p_2-k}\\\frac{(k+1)(p_2-k)+2p_1k}{kp_1t}J_{k\times p_1-k}^T& \frac{2}{p_1}(J_{p_1-k}-I_{p_1-k})& \frac{(p_1+p_2)(k+1)}{kp_1p_2}J_{p_1-k\times p_2-k}\\\frac{(k+1)(p_1-k)+2kp_2}{kp_2t}J_{k\times p_2-k}^T& \frac{(p_1+p_2)(k+1)}{kp_1p_2}J_{p_1-k\times p_2-k}^T& \frac{2}{p_2}(J_{p_2-k}-I_{p_2-k})\end{bmatrix}.$$
\end{proof}

\begin{example}
The \textit{kite graph} $Kite_{p,2}$ is obtained by identifying a vertex in $K_{p}$ to a pendant vertex of a path graph with $2$ vertices, which can be viewed as $K_{p}\circ_{1} K_{2}.$ Let $v^{*}$ be the identified vertex of a vertex $v_{1}$ in $K_{p}$ and the vertex $u_{1}$ in $K_{2}.$ Then by Theorem \ref{rd} we have
    \begin{itemize}
        \item [(a)]for $v_{i}= v^{*},\; v_{j}\in V(K_{p}\setminus \{v^{*}\}),\; r_{ij}= \frac{2}{p}.$
        \item[(b)] for $v_{i},v_{j}\in V(K_{p}\setminus \{v^{*}\}),\; r_{ij}= \frac{2}{p}.$
        \item[(c)] for $v_{i}= v^{*}$ and $v_{j}= u_{2}\in V(K_{2}),\; r_{ij}= 1.$
        \item[(d)] for $v_{i}\in V(K_{p}\setminus \{v^{*}\})$ and $v_{j}= u_2\in V(K_{2}),\; r_{ij}= \frac{p+2}{p}.$
    \end{itemize}
    
\end{example}

The windmill graph $W_{n+1}^{t}$ is the graph obtained by taking  $t\geq 2$ copies of complete graph $K_{n+1}$, for 
$n \geq 1$, with a vertex in common. By the definition of coalescence of graphs one can easily write $W_{n+1}^{t}= K_{n+1}\underbrace{\circ_{1}\cdots\circ_{1}}_{\substack{t-times}} K_{n+1}.$
Next theorem gives the resistance distance in $W_{n+1}^{t}$.
\begin{proposition}\label{wind}
   For $n>1,$ the resistance distance of vertices in $W_{n+1}^{t}$ is given by,
   $$r_{ij}= \begin{cases}
       \frac{4}{n+1},&\text{if $v_i, v_j$ are in different blocks,} \\
       \frac{2}{n+1}, &\text{otherwise}.
       
   \end{cases}$$
\end{proposition}

\begin{proof}
    The Laplacian matrix of $W_{n+1}^{t}$ is
    $$L(W_{n+1}^{t})= \begin{bmatrix}
        tnI_1& -J_{1\times tn}\\-J_{tn\times 1}& I_t\otimes (n+1I_n-J_n)
    \end{bmatrix}.$$
    Then its $\{1\}$-inverse is 
\begin{equation*}
    \begin{aligned}
        L^{(1)}(W_{n+1}^{t})&= \begin{bmatrix}
    \frac{1}{tn}I_1&0\\0& (I_t\otimes (n+1I_n-J_n)-\frac{1}{tn}J_{tn})^{\#}
\end{bmatrix}\\
&=\begin{bmatrix}
    \frac{1}{tn}I_1&0\\0& I_t\otimes \frac{1}{n+1}(I_n+J_n)-\frac{1}{tn}J_{tn}
\end{bmatrix}.
    \end{aligned}
\end{equation*}

    Now by the definition of resistance distance we get the required result.
\end{proof}

From Proposition \ref{wind}, we get the following corollaries.

\begin{corollary}
The Kirchhoff index of $W_{n+1}^{t}$ is
$$Kf(W_{n+1}^{t})= \frac{2n^2t^2-n^2t+nt}{n+1}.$$
\end{corollary}

\begin{corollary}
The Kemeny's constant of $W_{n+1}^{t}$ is
$$\kappa(W_{n+1}^{t})= \frac{n^2(2t-1)}{n+1}.$$
\end{corollary}

A $3$-rose graph is a graph consisting of three cycles intersecting in
a common vertex. Let $\mathcal{R}(r, s, t)$ denote the $3$-rose graph on $n = r + s + t - 2$ vertices, that is, the graph consisting of
three cycles $C_a$, $C_b$ and $C_c$ intersecting in a common vertex. Next corollary directly follows from Proposition \ref{wind}.
\begin{corollary}
    The resistance distance in $\mathcal{R}(3,3,3)$ is given by,
    $$r_{ij}= \begin{cases}
        \frac{4}{3}, &\text{if $v_i, v_j$ are in different blocks,} \\
       \frac{2}{3}, &\text{otherwise}.
    \end{cases}$$
\end{corollary}

\begin{theorem}
    Let $G_n$ be a graph of order $n$. For $p\geq k$, let $T$ be the collection of vertices in $K_{p}\circ_{k}(G_{n}\vee K_k)$, which are the identified vertices of $K_k$ and some verices in $K_{p}.$  Then,
  \begin{itemize}
    \item [(i)] for any $v_i, v_j\in T$,
    $r_{ij}= \frac{2}{p+n}.$
    \item [(ii)] for $v_i\in T, v_j\in V(K_{p}\setminus T)$,
    $r_{ij}= \frac{k(2p+n)+n}{kp(p+n)}.$
    \item [(iii)] for $v_i\in T, v_j\in V(G_{n})$,
    $r_{ij}= \frac{k-1}{k(p+n)}+(L(G_n)+kI_n)^{-1}_{jj}.$
    \item [(iv)] for any $v_i, v_j\in V(K_{p}\setminus T), r_{ij}= \frac{2}{p}.$
    \item [(v)] for $v_i\in V(K_{p}\setminus T), v_j\in V(G_{n}), r_{ij}= \frac{k+1}{kp}+(L(G_n)+kI_n)^{-1}_{jj}.$
    \item [(vi)] for any $v_i, v_j\in V(G_{n}), r_{ij}= (L(G_n)+kI_n)^{-1}_{ii}+(L(G_n)+kI_n)^{-1}_{jj}-2(L(G_n)+kI_n)^{-1}_{ij}.$
\end{itemize}
   
\end{theorem}
\begin{proof}
The Laplacian matrix of $K_p\circ_{k}(G_{n}\vee K_k)$ is given by,
$$L(K_p\circ_{k}(G_{n}\vee K_k))= \begin{bmatrix}(p+n)I_k-J_k& -J_{k\times p-k}& -J_{k\times n}\\-J_{k\times p-k}^T& p I_{p-k}-J_{p-k}& 0\\-J_{k\times n}^T& 0& L(G_n)+kI_n\end{bmatrix}.$$

Let $L_1= \begin{bmatrix}(p+n)I_k-J_k& -J_{k\times p-k}\\-J_{p-k\times k}& pI_{p-k}-J_{p-k}\end{bmatrix}, L_2= \begin{bmatrix}-J_{k\times n}\\0\end{bmatrix}$, and $L_3= L(G_n)+kI_n.$\\
Then by Lemmas \ref{P} and \ref{rs} we get,\\
$$L_{1}^{-1}= \begin{bmatrix}\frac{1}{p+n}(I_k+\frac{p}{nk}J_k)& \frac{1}{nk}J_{k\times p-k}\\ \frac{1}{nk}J_{p-k\times k} & \frac{1}{p}(I_{p-k}+\frac{p+n}{nk}J_{p-k})\end{bmatrix}.$$
Now let $S= L_3-L_{2}^{T}L_{1}^{-1}L_2$, then 
    $S= L(G_n)+kI_n-\frac{k}{n}J_n.$
   
From Lemma \ref{dag}, $S^\#= (L(G_n)+kI_n)^{-1}-\frac{k}{n}J_n.$\\
Therefore,
$$L^{(1)}(K_p\circ_{k}(G_{n}\vee K_k))= \begin{bmatrix}\frac{1}{p+n}(I_k+\frac{p}{nk}J_k)& \frac{1}{nk}J_{k\times p-k}& 0\\\frac{1}{nk}J_{p-k\times k} & \frac{1}{p}(I_{p-k}+\frac{p+n}{nk}J_{p-k})& 0\\0& 0& (L(G_n)+kI_n)^{-1}-\frac{1}{kn}J_n\end{bmatrix}.$$
By applying the definition of resistance distance, we obtain the required result.
\end{proof}
\begin{theorem}
     If $G_n$ is a graph of order $n$, then the resistance distance of the vertices in $K_{1, p-1}\circ_{1}(G_{n}\vee K_1)$ is given by,
     \begin{itemize}
\item[(i)] for $v_i=u^{*}, v_j\in V(K_{1,p-1}\setminus \{u^{*}\})$, $r_{ij}= 1,$
    \item [(ii)] for $v_i= u^{*}$, $v_j\in V(G_n)$, $r_{ij}= (L(G_{n})+I_n)^{-1}_{jj},$
    \item [(iii)]for $v_i, v_j\in V(K_{1,p-1}\setminus \{u^{*}\}), r_{ij}= 2,$
    \item [(iv)] for $v_i\in V(K_{1,p-1}\setminus \{u^{*}\}), v_j\in V(G_n)$, $r_{ij}= 1+(L(G_n)+I_n)^{-1}_{jj},$
    \item [(v)]for $v_i, v_j\in V(G_n), r_{ij}= (L(G_n)+I_n)^{-1}_{ii}+(L(G_n)+I_n)^{-1}_{jj}-2(L(G_n)+I_n)^{-1}_{ij},$
\end{itemize}
where $u^{*}$ is the identified vertex of a vertex in $K_1$ and a vertex in $K_{1,p-1}$ (center).
\end{theorem}
\begin{proof}
The Laplacian matrix of $K_{1, p-1}\circ_{1}(G_{n}\vee K_1)$ is given by,
$$L(K_{1, p-1}\circ_{1}(G_{n}\vee K_1))= \begin{bmatrix}L_1&L_2\\L_2^T&L_3\end{bmatrix},$$
where $L_1= \begin{bmatrix}(p+n-1)I_1&-J_{1\times p-1}\\-J_{p-1\times 1}&I_{p-1}\end{bmatrix},L_2= \begin{bmatrix}
    -J_{1\times n}\\0
\end{bmatrix}$, and $L_3= L(G_n)+I_n.$
\\
Then by Lemmas \ref{P} and \ref{rs} we get,\\
$$L_{1}^{-1}= \begin{bmatrix}\frac{1}{n}I_1& \frac{1}{n}J_{1\times p-1}\\ \frac{1}{n}J_{p-1\times 1} & I_p+\frac{1}{n}J_p\end{bmatrix}.$$
Now let $S= L_3-L_{2}^{T}L_{1}^{-1}L_2$, then $S= L(G_n)+I_n-\frac{1}{n}J_n.$\\
From Lemma \ref{dag}, $S^\#= (L(G_n)+I_n)^{-1}-\frac{1}{n}J_n.$\\
Therefore,
$$L^{(1)}(K_{1, p-1}\circ_{1}(G_{n}\vee K_1))= \begin{bmatrix}\frac{1}{n}I_1& \frac{1}{n}J_{1\times p-1}& 0\\\frac{1}{n}J_{p-1\times 1} & I_{p-1}+\frac{1}{n}J_{p-1}& 0\\0& 0& (L(G_n)+I_n)^{-1}-\frac{1}{n}J_n\end{bmatrix}.$$
Now by the definition of resistance distance we get the required result.
\end{proof}

\begin{theorem}
    The resistance distance matrix of $K_{p, q}\circ_{1}K_{1,n}$ is given by,
  $$R(K_{p, q}\circ_{1}K_{1,n})= \begin{bmatrix}
      0 &\frac{2}{q}J_{1\times p-1}& \frac{p+q-1}{pq}J_{1\times q}&J_{1\times n}\\ \frac{2}{q}J_{p-1\times 1}& \frac{2}{q}(J_{p-1}-I_{p-1})&\frac{p+q-1}{pq}J_{p-1\times q}&\frac{q+2}{2}J_{p-1\times n}\\ \frac{p+q-1}{pq}J_{q\times 1}&\frac{p+q-1}{pq}J_{q\times p-1}&\frac{2}{p}(J_q-I_q)&\frac{q(p+1)+(p-1)}{pq}J_{q\times n}\\J_{n\times 1}&\frac{q+2}{q}J_{n\times p-1}& \frac{q(p+1)+(p-1)}{pq}J_{n\times q}& 2(J_n-I_n)
  \end{bmatrix}.$$
    
\end{theorem}
\begin{proof}
    The Laplacian matrix of $K_{p, q}\circ_{1}K_{1,n}$ is given by,
$$L(K_{p, q}\circ_{1}K_{1,n})= \begin{bmatrix}L_1&L_2\\L_2^T&L_3\end{bmatrix},$$
where $L_1= \begin{bmatrix}(q+n)I_1&0&-J_{1\times q}\\0&qI_{p-1}&-J_{p-1\times q}\\-J_{q\times 1}&-J_{q\times p-1}&pI_q\end{bmatrix},L_2= \begin{bmatrix}
    -J_{1\times n}\\0\\0\end{bmatrix}$, and $L_3= I_{n}.$
\\
Then by Lemmas \ref{P} and \ref{rs} we get,\\
$$L_{1}^{-1}= \begin{bmatrix}\frac{1}{n}I_1&J_{1\times p-1}& \frac{1}{n}J_{1\times q}\\ \frac{1}{n}J_{p-1\times 1}&\frac{1}{q}I_{p-1}+\frac{n+q}{nq}J_{p-1}&\frac{q+n}{nq}J_{p-1\times q}\\\frac{1}{n}J_{q\times 1}&\frac{q+n}{nq}J_{q\times p-1}& \frac{1}{p}I_q+\frac{p(q+n)-n}{npq}J_q\end{bmatrix}.$$
Consider \begin{equation*}
    \begin{aligned}
        S&= L_3-L_{2}^{T}L_{1}^{-1}L_2\\
        &= I_{n}-\frac{1}{n}J_{n}
    \end{aligned}
\end{equation*}
then, $S^\#= I_{n}-\frac{1}{n}J_{n}.$\\
Therefore,
$$\scriptsize{L^{(1)}(K_{p, q}\circ_{1}K_{1,n})= \begin{bmatrix}\frac{1}{n}I_1&J_{1\times p-1}& \frac{1}{n}J_{1\times q}&0\\ \frac{1}{n}J_{p-1\times 1}&\frac{1}{q}I_{p-1}+\frac{n+q}{nq}J_{p-1}&\frac{q+n}{nq}J_{p-1\times q}&0\\\frac{1}{n}J_{q\times 1}&\frac{q+n}{nq}J_{q\times p-1}& \frac{1}{p}I_q+\frac{p(q+n)-n}{npq}J_q&0\\0&0&0&I_n-\frac{1}{n}J_n\end{bmatrix}}.$$
Now by the definition of resistance distance we get the required result.
\end{proof}
\begin{theorem}\label{kpq}
    The resistance distance matrix of $K_{p, q}\circ_{1}K_n$ is given by,
  $$R(K_{p, q}\circ_{1}K_n)= \begin{bmatrix}
      0&\frac{2}{q}J_{1\times p-1}&\frac{p+q-1}{pq}J_{1\times q}&\frac{2}{n}J_{1\times n}\\ \frac{2}{q}J_{p-1\times 1}&\frac{2}{q}(J_{p-1}-I_{p-1})&\frac{p+q-1}{pq}J_{p-1\times q}&\frac{2(q+n)}{qn}J_{p-1\times n}\\ \frac{p+q-1}{pq}J_{q\times 1}&\frac{p+q-1}{pq}J_{q\times p-1}&\frac{2}{p}(J_q-I_q)&\frac{q(n+2p)+n(p-1)}{npq}J_{q\times n}\\ \frac{2}{n}J_{n\times 1}& \frac{2(q+n)}{qn}J_{n\times p-1}&\frac{q(n+2p)+n(p-1)}{npq}J_{n\times q}&\frac{2}{n}(J_{n-1}-I_{n-1})
  \end{bmatrix}.$$
   
\end{theorem}
\begin{proof}
    The Laplacian matrix of $K_{p, q}\circ_{1}K_n$ is given by,
$$L(K_{p, q}\circ_{1}K_n)= \begin{bmatrix}L_1&L_2\\L_2^T&L_3\end{bmatrix},$$
where $L_1= \begin{bmatrix}(q+n-1)I_1&0&-J_{1\times q}\\0&qI_{p-1}&-J_{p-1\times q}\\-J_{q\times 1}& -J_{q\times p-1}& pI_q\end{bmatrix},L_2= \begin{bmatrix}
    -J_{1\times n-1}\\0\\0
\end{bmatrix}$, and $L_3= nI_{n-1}-J_{n-1}.$
\\
Then by Lemmas \ref{P} and \ref{rs} we get,\\
$$L_{1}^{-1}= \begin{bmatrix}\frac{1}{n-1}I_1&\frac{1}{n-1}J_{1\times p-1}& \frac{1}{n-1}J_{1\times q}\\ \frac{1}{n-1}&\frac{1}{q}(I_{p-1}+\frac{q+n-1}{n-1}J_{p-1})&\frac{q+n-1}{q(n-1)}J_{p-1\times q}\\ \frac{1}{n-1}J_{q\times 1}&\frac{q+n-1}{q(n-1)}J_{q\times p-1}& \frac{1}{p}(I_q+\frac{p(q+n-1)-(n-1)}{q(n-1)}J_q)\end{bmatrix}.$$
Consider \begin{equation*}
    \begin{aligned}
        S&= L_3-L_{2}^{T}L_{1}^{-1}L_2\\
        &= nI_{n-1}-\frac{n}{n-1}J_{n-1}
    \end{aligned}
\end{equation*}
then, $S^\#= \frac{1}{n}I_{n-1}-\frac{1}{n(n-1)}J_{n-1}.$\\
Therefore,
$$\scriptsize{L^{(1)}(K_{p, q}\circ_{1}K_n)= \begin{bmatrix}\frac{1}{n-1}I_1&\frac{1}{n-1}J_{1\times p-1}& \frac{1}{n-1}J_{1\times q}&0\\ \frac{1}{n-1}&\frac{1}{q}(I_{p-1}+\frac{q+n-1}{n-1}J_{p-1})&\frac{q+n-1}{q(n-1)}J_{p-1\times q}&0\\ \frac{1}{n-1}J_{q\times 1}&\frac{q+n-1}{q(n-1)}J_{q\times p-1}& \frac{1}{p}(I_q+\frac{p(q+n-1)-(n-1)}{q(n-1)}J_q)&0\\0&0&0&\frac{1}{n}(I_{n-1}-\frac{1}{n-1}J_{n-1})\end{bmatrix}}.$$
Now by the definition of resistance distance we get the required result.
\end{proof}

The \textit{pineapple graph} $K_{p}^{q}$ is the coalescence of the complete graph $K_p$ (at any vertex) with
the star $K_{1,q}$ at the vertex of degree $q.$ It has $n = p + q$ vertices and $^pC_{2}+q$ edges. \begin{figure}[htbp]
\centering
\begin{tikzpicture}[scale=0.4]
 \filldraw[ draw = black] (-15.8,1.4) circle (0.2 cm);
 \filldraw[ draw = black] (-8.3,1.4) circle (0.2 cm);
 \filldraw[ draw = black] (-15.8,-1.5) circle (0.2 cm);
 \filldraw[ draw = black] (-8.3,-1.5) circle (0.2 cm);
       \filldraw[ draw = black] (-12,4) circle (0.2 cm);
 \filldraw[ draw = black] (-12,-4) circle (0.2 cm);
\filldraw[ draw = black] (-12,6) circle (0.2 cm);
 \filldraw[ draw = black] (-14,6) circle (0.2 cm);
\filldraw[ draw = black] (-10,6) circle (0.2 cm);
 \filldraw[ draw = black] (-8,6) circle (0.2 cm);
  \filldraw[ draw = black] (-16,6) circle (0.2 cm);
\draw (-16,6)--(-12,4)--(-8,6);
\draw (-10,6)--(-12,4)--(-14,6);
\draw (-12,4)--(-12,6);
\draw (-12,4)--(-8.3,-1.5);
\draw (-15.8,1.4)--(-8.3,-1.5)--(-15.8,-1.5)--(-15.8,1.4)--(-12,4)--(-8.3,1.4)--(-15.8,1.4)--(-12,-4);
\draw (-15.8,-1.5)--(-12,4)--(-12,-4)--(-15.8,1.4)--(-15.8,-1.5)--(-12,-4)--(-8.3,1.4)--(-8.3,-1.5)--(-12,-4)--(-15.8,-1.5)--(-8.3,1.4)--(-15.8,1.4);
 \end{tikzpicture}
\caption{$K_{6}^{5}$.}
 \end{figure}
 
 Using Theorem \ref{kpq}, we get the following proposition.
\begin{proposition}\label{pin}
The resistance distance in a pineapple graph $K_{p}^{q}$ is given by
\begin{itemize}
\item[(i)] for $v_i=v^{**}, v_j\in V(K_{p}\setminus \{v^{**}\})$, $r_{ij}= \frac{2}{p},$
    \item [(ii)] for $v_i= v^{**}$, $v_j\in V(K_{1,q}\setminus \{v^{**}\})$, $r_{ij}= 1,$
    \item [(iii)] for $v_i, v_j\in V(K_{p}\setminus \{v^{**}\}), r_{ij}= \frac{2}{p},$
    \item [(iv)] for $v_i\in V(K_{p}\setminus \{v^{**}\}), v_j\in V(K_{1,q}\setminus \{v^{**}\})$, $r_{ij}= \frac{p+2}{p},$
    \item [(v)] for $v_i, v_j\in V(K_{1,q}\setminus \{v^{**}\}), r_{ij}= 2,$
\end{itemize}
where $v^{**}$ is the identified vertex of a vertex $u$ in $K_p$ and the vertex $v$ of degree $q$ in $K_{1,q}$.
\end{proposition}

From Proposition \ref{pin}, we get the following corollaries.

\begin{corollary}
The Kirchhoff index of $K_{p}^{q}$ is
$$Kf(K_{p}^{q})= q(p+q+3)+p+1-\frac{2}{p}(q+1).$$
\end{corollary}

\begin{corollary}
The Kemeny's constant of $K_{p}^{q}$ is
$$\kappa(K_{p}^{q})= \frac{p^4-p^3+p^3q+3p^2q+2pq^2-3p^2-7pq+7p+4q-2}{p(p-1)+2q}.$$
\end{corollary}

The \textit{dandelion graph} $D(n,l)$ on $n$ vertices is the coalescence of the star graph $K_{1,n-l}$ (at the center) with
the path $P_{l}$ at any pendant vertex. 
\begin{figure}[htbp]
\centering
\begin{tikzpicture}[scale=0.4]

     \filldraw[ draw = black] (1.5,3.7) circle (0.2 cm);
      \filldraw[ draw = black] (3.7,1.6) circle (0.2 cm);
       \filldraw[ draw = black] (3.7,-1.5) circle (0.2 cm);
        \filldraw[ draw = black] (1.5,-3.7) circle (0.2 cm);
         \filldraw[ draw = black] (-1.5,-3.7) circle (0.2 cm);
          \filldraw[ draw = black] (-3.7,-1.4) circle (0.2 cm);
           \filldraw[ draw = black] (-3.8,1.4) circle (0.2 cm);
            \filldraw[ draw = black] (-1.6,3.7) circle (0.2 cm);
  \filldraw[ draw = black] (0,4) circle (0.2 cm);
   \filldraw[ draw = black] (0,-4) circle (0.2 cm);
    \filldraw[ draw = black] (4,0) circle (0.2 cm);
     \filldraw[ draw = black] (-4,0) circle (0.2 cm);
     \filldraw[ draw = black] (0,0) circle (0.2 cm);   

               \filldraw[ draw = black] (2.8,2.8) circle (0.2 cm);  
     \filldraw[ draw = black] (-2.8,-2.8) circle (0.2 cm);  
        \filldraw[ draw = black] (2.8,-2.8) circle (0.2 cm);
         \filldraw[ draw = black] (-2.8,2.8) circle (0.2 cm);
          \filldraw[ draw = black] (8,0) circle (0.2 cm);
           \filldraw[ draw = black] (12,0) circle (0.2 cm);
           \draw (4,0)--(8,0)--(12,0);
         \draw (2.8,2.8)--(0,0)--(-2.8,-2.8);
         \draw (2.8,-2.8)--(0,0)--(-2.8,2.8);
\draw (0,4)--(0,0)--(0,-4);
\draw (-1.6,3.7)--(0,0)--(-3.8,1.4);          
\draw (4,0)--(0,0)--(-4,0);
\draw (-3.7,-1.4)--(0,0)--(-1.5,-3.7);
\draw (1.5,-3.7)--(0,0)--(3.7,-1.5);
\draw (3.7,1.6)--(0,0)--(1.5,3.7);
       
\end{tikzpicture}
\caption{$D(19,4)$.}
 \end{figure}
 
 The following theorem describes the resistance distance matrix of a dandelion graph $D(n,l).$

\begin{theorem}\label{dand}
The resistance distance matrix of a dandelion graph $D(n,l)$  on $n$ vertices is given by\\
$$R(D(n,l))= \begin{bmatrix}
    0& 1& 2& \cdots & l-1& 1&& \cdots & 1\\
1& 0 & 1 &\cdots& l-2& 2& &\cdots &2 \\
2&1&0&\cdots&l-3&3&&\cdots&3\\
\vdots & \vdots &\ddots&\ddots & \vdots & \vdots & &\vdots & \vdots  \\
l-1&l-2& \cdots&1&0&l&&\cdots &l\\
1&2&&\cdots& l&0& 2&\cdots &2\\
1&2&&\cdots& l&2&0&\cdots &2\\
\vdots & \vdots& &\ddots & \vdots & \vdots & \vdots & \ddots & \vdots \\
1&2&&\cdots&l&2&2&\cdots&0

\end{bmatrix}.$$
\end{theorem}
\begin{proof}
By a proper labelling of vertices in $D(n,l)= K_{1,n-l}\circ_{1}P_{l}$, we can write its Laplacian matrix as\\
$$L(D(n,l))= \begin{bmatrix}L_{1}& L_{2}\\
L_{2}^{T}&L_{3}.\end{bmatrix},$$\\
where $L_1= \begin{bmatrix}
    n+1-l&-1&0&\cdots&0\\
    -1&2&-1&\cdots&0\\
    \vdots&\cdots&\ddots&\vdots&\vdots\\
    0&\cdots&-1&2&-1\\
    0&\cdots&0&-1&1
    
\end{bmatrix}_{l\times l}$, $L_2= \begin{bmatrix}
    -J_{1\times n-l}\\
    0_{l-1\times n-l}
\end{bmatrix}$ and $L_3= I_{n-l}.$\\
Then by applying Lemma \ref{L} we get,\\
$$L^{(1)}= \begin{bmatrix}L_{1}^{-1}&0\\
0& I_{n-l}-J_{n-l}\end{bmatrix},$$
where $L_{1}^{-1}= \begin{bmatrix}
    \frac{1}{n-l}&\frac{1}{n-l}&\cdots&\frac{1}{n-l}\\
    \frac{1}{n-l}&\frac{(n-l)+1}{n-l}&\cdots&\frac{(n-l)+1}{n-l}\\
    \vdots&\vdots&\ddots&\vdots\\
    \frac{1}{n-l}&\frac{(n-l)+1}{n-l}&\cdots&\frac{(l-1)(n-l)+1}{n-l}
\end{bmatrix}.$
\\
\\
Now by the definition of resistance distance we get the required result.
\end{proof}

From Theorem \ref{dand}, we get the following corollaries.

\begin{corollary}
The Kirchhoff index of $D(n,l)$ is
$$Kf(D(n,l))= \frac{l^2(3n-2l+2)+l(5-9n)+6(n^2-1)}{6}.$$
\end{corollary}

\begin{corollary}
The Kemeny's constant of $D(n,l)$ is
$$\kappa(D(n,l))= \frac{(n+1)(2l^2-1)+2n(n-3l)}{2(n-1)}+\frac{l(5-2l^2)}{3(n-1)}.$$
\end{corollary}

The following theorems describes various graph parameters of $K_{p_1}\circ_{k}K_{p_2}.$
\begin{proposition}\label{kir}
The Kirchhoff index of $K_{p_1}\circ_{k}K_{p_2}$ is given by
\begin{equation*}
 \begin{aligned}
   \mathcal{K}f(K_{p_1}\circ_{k}K_{p_2})&= \frac{1}{kp_1p_2t}\bigl((p_1-k)(p_2-k)(k+1)\left(p_1(k+t)+p_2+\frac{p_2k}{k+1}+\frac{p_2k(p_1k+p_1t-t)}{(p_2-k)(k+1)}\right)\\&+kp_1(p_2-k)\left(p_2(t+2k)-t(k+1)+kp_2(k-1)\right)\bigr).
 \end{aligned}
 \end{equation*}

\end{proposition}
\begin{proof}
 By definition $\mathcal{K}f(G_n)= \sum_{i<j} r_{ij}(G_n).$ Then,
\begin{equation*}
    \begin{aligned}
      \mathcal{K}f(K_{p_1}\circ_{k}K_{p_2})&= \sum_{v_i, v_j\in T}\frac{2}{t}+ \sum_{v_i\in T, v_j\in V(K_{p_1}\setminus T)}\frac{k(p_1+t)+(p_2-k)}{p_1kt}\\&\;\;+\sum_{v_i\in T, v_j\in V(K_{p_2}\setminus T)}\frac{(k+1)(p_1-k)+2p_2k}{p_2kt}+ \sum_{v_i, v_j\in V(K_{p_1}\setminus T)}\frac{2}{p_1}\\&\;\;+\sum_{v_i\in V(K_{p_1}\setminus T), v_j\in V(K_{p_2}\setminus T)}\frac{(p_1+p_2)(k+1)}{kp_1p_2}+\sum_{v_i, v_j\in V(K_{p_2}\setminus T)}\frac{2}{p_2}.\\
      &= \frac{1}{kp_1p_2t}\bigl((p_1-k)(p_2-k)(k+1)\left(p_1(k+t)+p_2+\frac{p_2k}{k+1}+\frac{p_2k(p_1k+p_1t-t)}{(p_2-k)(k+1)}\right)\\&\;\;+kp_1(p_2-k)\left(p_2(t+2k)-t(k+1)+kp_2(k-1)\right)\bigr).
    \end{aligned}
\end{equation*}
\end{proof}
\begin{proposition}\label{kappa}
The Kemeny's constant of $K_{p_1}\circ_{k}K_{p_2}$ is given by
{\scriptsize{\begin{equation*}
 \begin{aligned}
   \kappa(K_{p_1}\circ_{k}K_{p_2})= &\frac{1}{2mp_1p_2}\left(\frac{t-1}{t}(p_1(t-1)(p_2k(k-1)+(p_2-1)(p_2-k)((k+1)(p_1-k)+2p_2k))\right.\\&\left.+p_2(p_1-1)(p_1-k)(2p_1k+(p_2-k)(k+1)))+p_1(p_2-k)(p_2-1)^2(p_2-k-1)\right.\\&\left.+\frac{(p_1-k)(p_1-1)(p_2k(p_1-k-1)(p_1-1)+(p_2-k)(p_2-1)(p_1+p_2)(k+1))}{k}\right).
 \end{aligned}
 \end{equation*}}}
 
\end{proposition}

\begin{proposition}
 The additive degree-Kirchhoff index of $K_{p_1}\circ_{k}K_{p_2}$ is given by
\scriptsize{\begin{equation*}
 \begin{aligned}
   R^{+}(K_{p_1}\circ_{k}K_{p_2})= &\frac{2p_1k(k-1)(t-1)+k(p_1-k)(2t-2)\left((2t)+p_2-k\right)}{p_1t}+ \frac{(p_1-k)(p_2-k)(p_1+p_2-2)(p_1+p_2)(k+1)}{kp_1p_2}\\&\;\;+\frac{(p_2-k)(p_1+2p_2-k-2)\left((k+1)(p_1-k)+2p_2k\right)}{p_2t}\\&\;\;
   +\frac{2p_2(p_1-k)(p_1-k-1)(p_1-1)+2p_1(p_2-k)(p_2-k-1)(p_2-1)}{p_1p_2}.
 \end{aligned}
 \end{equation*}}   
\end{proposition}

\begin{proposition}
    The multiplicative degree-Kirchhoff index of $K_{p_1}\circ_{k}K_{p_2}$ is given by
{\scriptsize{\begin{equation*}
 \begin{aligned}
   R^{*}(K_{p_1}\circ_{k}K_{p_2})= &\frac{1}{p_1p_2}\left(\frac{t-1}{t}(p_1(t-1)(p_2k(k-1)+(p_2-1)(p_2-k)((k+1)(p_1-k)+2p_2k))\right.\\&\left.+p_2(p_1-1)(p_1-k)(2p_1k+(p_2-k)(k+1)))+p_1(p_2-k)(p_2-1)^2(p_2-k-1)\right.\\&\left.+\frac{(p_1-k)(p_1-1)(p_2k(p_1-k-1)(p_1-1)+(p_2-k)(p_2-1)(p_1+p_2)(k+1))}{k}\right).
 \end{aligned}
 \end{equation*}}} 
 
\end{proposition}

\begin{proposition}
    The mixed degree-Kirchhoff index of $K_{p_1}\circ_{k}K_{p_2}$ is given by
\scriptsize{\begin{equation*}
 \begin{aligned}
   \hat{R}(K_{p_1}\circ_{k}K_{p_2})= &\frac{2k(k-1)}{t}+\frac{(p_1-k)((t-1)^2+(p_1-1)^2)(k(2t)+(p_2-k))}{p_1(p_1-1)t(t-1)}\\
   &+\frac{(p_2-k)((t-1)^2+(p_2-1)^2)((k+1)(p_1-k)+2p_2k)}{p_2(p_2-1)t(t-1)}\\
   &+\frac{2(p_1-k)(p_1-k-1)}{p_1}+\frac{2(p_2-k)(p_2-k-1)}{p_2}\\
   &+\frac{(p_1-k)(p_2-k)((p_1-1)^2+(p_2-1)^2)(p_1+p_2)(k+1)}{kp_1p_2(p_1-1)(p_2-1)}.
 \end{aligned}
 \end{equation*}}
\end{proposition}

In general, it is difficult to find the resistance energy of graphs. The following table gives the resistance energy of $K_{p_1}\circ_{k}K_{p_2}$, for different values of $p_1, p_2$ and $k.$\\
   \begin{minipage}{0.32\textwidth}
       \begin{table}[H]
\begin{tabular}{@{}lllll@{}}
\toprule

No. &$p_1$& $p_2$&$k$& $RE(K_{p_1}\circ_{k}K_{p_2})$\\
 
\midrule

	$1$&3&2&1&6.21\\
$2$& 4&2&1&6.92\\
 $3$&5&2&1&7.79\\
 $4$&6&2&1&7.91\\
 $5$&2&3&1&6.23\\
 $6$&3&3&1&7.29\\
 $7$&4&3&1&7.91\\
 8&5&3&1&8.37\\
 9&6&3&1&8.76\\
 10&2&4&1&6.92\\
 11&3&4&1&7.9\\
 12&4&4&1&8.44\\
 13&5&4&1&8.82\\
 14&6&4&1&9.14\\
 15&5&5&1&9.13\\
 16&6&5&1&9.4\\
 17&7&5&1&9.6\\
 18&2&2&2&4\\
 19&3&2&2&4.36\\
 20&4&2&2&4.45\\
 21&5&2&2&4.48\\
 22&6&2&2&4.48\\
23&2&3&2&4.4\\
	
\bottomrule

\end{tabular}

\end{table}

  \strut \end{minipage}
  \hfill\allowbreak
  \begin{minipage}{0.62\textwidth}
       \begin{table}[H]
\begin{tabular}{@{}lllll@{}}
\toprule

 No. &$p_1$& $p_2$&$k$& $RE(K_{p_1}\circ_{k}K_{p_2})$\\
 
\midrule
 
 24&3&3&2&4.92\\
 25&4&3&2&5.63\\
 26&5&3&2&5.42\\
 27&6&3&2&5.59\\
 28&2&4&2&4.56\\
 29&3&4&2&5.23\\
 30&4&4&2&5.64\\
 31&5&4&2&5.92\\
 32&6&4&2&6.17\\
 33&5&5&2&6.27\\
 34&6&5&2&6.52\\
 35&7&5&2&6.73\\
 36&4&3&3&8.14\\
 37&5&3&3&8.18\\
 38&6&3&3&8.19\\
 39&4&4&3&8.25\\
 40&5&4&3&8.3\\
 41&6&4&3&8.43\\
 42&5&5&3&8.36\\
 43&6&5&3&8.39\\
 44&7&5&3&8.41\\
 45&6&5&4&2.92\\
 46&7&5&4&2.77\\
	
\bottomrule

\end{tabular}

\end{table}
\strut \end{minipage}

\section{Conclusion}
    This article explores the concept of resistance distance in the $k$-coalescence of complete graphs. These results enable us to determine several graph parameters, including Kemeny's constant, Kirchhoff index, additive degree-Kirchhoff index, multiplicative degree-Kirchhoff index, and mixed degree-Kirchhoff index of $k$-coalescence of complete graphs. In addition, the resistance distance in the $k$-coalescence of a complete graph with particular graphs are obtained. Furthermore, the article applies the findings to determine the resistance distance of specific graphs like the vertex coalescence of a complete bipartite graph with a complete graph, a complete bipartite graph with a star graph, the windmill graphs, the pineapple graph, etc. 

\section{Declarations}
 On behalf of all authors, the corresponding author states that there is no conflict of interest.

 \bibliography{rd}

\begin{thebibliography}{10}

\bibitem{bapat1999resistance}
R~B Bapat.
\newblock Resistance distance in graphs.
\newblock {\em Mathematics Student-India}, 68(1-4):87--98, 1999.

\bibitem{bapat2010graphs}
R~B Bapat.
\newblock {\em Graphs and matrices}, volume~27.
\newblock Springer, 2010.

\bibitem{bapat2010resistance}
R~B Bapat and S~Gupta.
\newblock Resistance distance in wheels and fans.
\newblock {\em Indian Journal of Pure and Applied Mathematics}, 41(1):1--13,
  2010.

\bibitem{bapat2003simple}
R~B Bapat, I~Gutman, and W~Xiao.
\newblock A simple method for computing resistance distance.
\newblock {\em Zeitschrift f{\"u}r Naturforschung A}, 58(9-10):494--498, 2003.

\bibitem{ben2003generalized}
A~Ben-Israel and T~N~E Greville.
\newblock {\em Generalized inverses: theory and applications}, volume~15.
\newblock Springer Science \& Business Media, 2003.

\bibitem{bianchi2016upper}
M~Bianchi, A~Cornaro, Jos{\'e}~Luis P, and Anna T.
\newblock Upper and lower bounds for the mixed degree-{K}irchhoff index.
\newblock {\em Filomat}, 30(9):2351--2358, 2016.

\bibitem{bonchev1994molecular}
D~Bonchev, A~T Balaban, X~Liu, and D~J Klein.
\newblock Molecular cyclicity and centricity of polycyclic graphs. {I}.
  cyclicity based on resistance distances or reciprocal distances.
\newblock {\em International journal of quantum chemistry}, 50(1):1--20, 1994.

\bibitem{subdivision}
C~Bu, Bo~Yan, X~Zhou, and J~Zhou.
\newblock Resistance distance in subdivision-vertex join and subdivision-edge
  join of graphs.
\newblock {\em Linear Algebra and its Applications}, 458:454--462, 2014.

\bibitem{chen2007resistance}
H~Chen and F~Zhang.
\newblock Resistance distance and the normalized {L}aplacian spectrum.
\newblock {\em Discrete applied mathematics}, 155(5):654--661, 2007.

\bibitem{cvet}
D.~M. Cvetkovi{\'{c}}{,}~M. Doob and H.~Sachs.
\newblock {\em Spectra of graphs- Theory and applications, $3^{rd}$ edition}.
\newblock Johann Ambrosius Barth Verlag, Heidelberg, Leipzig, 1995.

\bibitem{evans2022algorithmic}
EJ~Evans and AE~Francis.
\newblock Algorithmic techniques for finding resistance distances on structured
  graphs.
\newblock {\em Discrete Applied Mathematics}, 320:387--407, 2022.

\bibitem{gutman2012degree}
I~Gutman, L~Feng, and G~Yu.
\newblock Degree resistance distance of unicyclic graphs.
\newblock {\em Transactions on Combinatorics}, 1(2):27--40, 2012.

\bibitem{kcoal}
Sudhir~R. J and Shrinath~L. P.
\newblock Spectra of $k$ coalescence of complete graphs.
\newblock {\em Asia Mathematika}, 5:113--118, 2021.

\bibitem{jahfar2020central}
T~K Jahfar and A~V Chithra.
\newblock Central vertex join and central edge join of two graphs.
\newblock {\em AIMS Mathematics}, 5(6):7214--7234, 2020.

\bibitem{kirkland2012group}
S~J Kirkland and M~Neumann.
\newblock {\em Group inverses of {M}-matrices and their applications}.
\newblock CRC Press, 2012.

\bibitem{klein1993resistance}
D~J Klein and M~Randi{\'c}.
\newblock Resistance distance.
\newblock {\em Journal of mathematical chemistry}, 12:81--95, 1993.

\bibitem{liu2015resistance}
X~Liu, J~Zhou, and C~Bu.
\newblock Resistance distance and {K}irchhoff index of {R}-vertex join and
  {R}-edge join of two graphs.
\newblock {\em Discrete Applied Mathematics}, 187:130--139, 2015.

\bibitem{wang2019resistance}
L~Sun, Z~Shang, and C~Bu.
\newblock Resistance distance and kirchhoff index of the {Q}-vertex (or edge)
  join graphs.
\newblock {\em Discrete Mathematics}, 344(8):112433, 2021.

\bibitem{zhang2007resistance}
H~Zhang and Y~Yang.
\newblock Resistance distance and {K}irchhoff index in circulant graphs.
\newblock {\em International Journal of Quantum Chemistry}, 107(2):330--339,
  2007.

\bibitem{zhou2014line}
J~Zhou, L~Sun, W~Wang, and C~Bu.
\newblock Line star sets for {L}aplacian eigenvalues.
\newblock {\em Linear Algebra and its Applications}, 440:164--176, 2014.

\end{thebibliography}
 \bibliographystyle{plain}
\end{document}